   \documentclass[twoside,reqno,11pt]{fcaa-var} %

\usepackage{graphicx}
\usepackage{epsfig}
\usepackage{amsthm}
\usepackage{amsmath}
\usepackage{latexsym}
\usepackage{amsfonts}
\usepackage{amssymb}

\newtheoremstyle{theorem}
  {15pt}          
  {15pt}  
  {\sl}  
  {\parindent}
  {\sc}  
  {. }   
  { }    
  {}     
\theoremstyle{theorem}
\newtheorem{lemma}{Lemma}[section]
\newtheorem{theorem}{Theorem}[section]

\newtheoremstyle{defi}
  {15pt}          
  {15pt}  
  {\rm}  
  {\parindent}     
  {\sc}  
  {. }    
  { }    
  {}     
\theoremstyle{defi}
\newtheorem{definition}{Definition}[section]
\newtheorem{remark}{Remark}[section]

 \def\proofname{\indent\rm P\ r\ o\ o\ f.}
 \def\proofend{\hfill$\Box$}
 \def\qed{\hfill$\Box$} 

 \usepackage{hyperref} 

 \def\theequation{\arabic{section}.\arabic{equation}}

 \def\themycopyrightfootnote{\vspace*{3pt}
 \copyright \, Year\,  Diogenes  Co., Sofia
   \par  \noindent pp. xxx--xxx, DOI: ......................
   \hfill  \vspace*{-36pt}
  }

 \title[SOLUTIONS OF THE MAIN BAUNDARY VALUE PROBLEMS \dots]{SOLUTIONS OF THE MAIN BAUNDARY VALUE PROBLEMS FOR THE TIME-FRACTIONAL TELEGRAPH EQUATION BY THE GREEN FUNCTION METHOD }
 \author[\normalsize M. O. Mamchuev]{\normalsize Murat O. Mamchuev}

\begin{document}

 \vbox to 2.5cm { \vfill }


 \bigskip \medskip

 \begin{abstract}


The inhomogenous time-fractional telegraph equation with Caputo derevatives with constant coefficients is considered.
For considered equation the general representation of regular solution in rectangular domain is obtained, and the fundamental solution is presented.
Using this representation and the properties of fundamental solution, the Cauchy problem and the basic problems in half-strip and rectangular domains are studied. 
For Cauchy problem the theorems of existence and uniqueness of solution are proved, and the explicit form of solution is constructed. 
The solutions of the investigated problems are constructed in terms of the appropriate Green functions, which are also constructed an explicit form.

 \medskip


 \smallskip


 \end{abstract}

 \maketitle

 \vspace*{-16pt}



\section{Introduction}\label{sec:1}

\setcounter{section}{1}
\setcounter{equation}{0}\setcounter{theorem}{0}


Consider the equation
\begin{equation}\label{mam1}
{\bf L} u(x,y)\equiv \partial_{0y}^{\alpha}u(x,y)+b\partial_{0y}^{\beta}u(x,y)
-\frac{\partial^2}{\partial x^2}u(x,y)+cu(x,y)=f(x,y),
\end{equation}
where  $\alpha=2\beta\in]0,2[,$
$b$ and $c$ are given real numbers, 
$f(x,y)$ is a given real function,
$\partial_{0y}^{\nu}$ is the Caputo fractional differentiation operator of order $\nu$ 
\cite [p. 11]{n}: 
$$
\partial_{ay}^{\nu}g(y)={\mathop{\rm sgn}}^n(y-a)D_{ay}^{\nu-n}g^{(n)}(y),
\quad n-1<\nu\leq n, \quad n\in {\mathbb N},
$$                                
$D_{ay}^{\nu}$ is the Riemann-Liouville fractional integro-differentiation operator
of order $\nu$  \cite [p. 9]{n}: 
$$D_{ay}^{\nu}g(y)=\frac{\mathop{\rm sgn}(y-a)}{\Gamma(-\nu)}
\int\limits_a^y \frac{g(s)ds}{|y-s|^{\nu+1}}, \quad  \nu< 0,$$
for $\nu \geq 0$ the operator $D_{ay}^{\nu}$ can be determined by 
recursive relation
$$D_{ay}^{\nu}g(y)=\mathop{\rm sgn}(y-a)\frac{d}{dy}D_{ay}^{\nu-1}g(y), \quad
\nu \geq 0,
$$
$\Gamma(z)$ is the gamma-function.
 
Equations of the form (\ref{mam1}) are used to describe anomalous diffusion processes observed 
in experiments related to blood circulation \cite{Cascaval-2002},
iterated Brownian motion, and telegraph processes with Brownian time
\cite{Orsinger-2003}, \cite{Orsinger-2004}.
The solution of Cauchy problem for the equation 
\begin{equation} \label{bib42}
\partial_{0t}^{\alpha}u(x,t)+b\partial_{0t}^{\beta}u(x,t)=
c^2\frac{\partial^2}{\partial x^2}u(x,t), 
\end{equation}
where  $0<\alpha=2\beta< 2,$
was obtained in \cite{Orsinger-2004} in terms of the inverse Fourier transform.
A closed-form expression for the solution was obtained for $\beta=1/2.$
Cauchy, signalling and first boundary value problems for the equation
(\ref{bib42}) were studied in \cite{Chen-2008}, \cite{Huang-2009} by separation of variables
and in \cite{Huang-2009} by the methods Laplace and Fourier integral transforms.  

In paper \cite{atan2007}  the generalized telegraph equation
(\ref{bib42}) was investigated for the case when  $0<\beta\leq \alpha< 2.$
The solutions of signalling and Cauchy problems and of a problem in bounded domain 
with homogenous boundary condition were obtained by using sepatation of variables
and Laplace trasform methods. For the case when $0<\beta,\alpha< 1$ a maximum principle
for solution of equation (\ref{bib42}) in bounded domain was proved. 
In paper \cite{Bazhlekova-2013} a nonlocal boundary value problem for
the equation  (\ref{bib42}) with $0<\beta\leq \alpha< 2$ was studied.
The solution is given as an expansion on the generalized eigenfunctions.

Note that the following 
fractional advection-dispersion equation
\begin{equation} \label{fade}
\partial_{0t}^{\alpha}u(x,t)-D\frac{\partial^2}{\partial x^2}u(x,t)
+\nu \frac{\partial}{\partial x}u(x,t) +\mu u(x,t)=0, 
\end{equation}
where $D>0,$ $\nu,$ and $\mu$ is the constants, by the change of unknown function,
can be reduced to equation (\ref{mam1}) with $b=0.$ 
By using Mellin and Laplace transforms,
and properties of H-functions, authors of paper \cite{Liu-2003}
derives the complete solution of equation (\ref{fade}).
In paper \cite{Huang-2005} for equation (\ref{fade}) the solutions of problems in half-space and bounded domains
were obtained in cases $\mu=0$ and $\nu=0$ respectively.

Note also that the boundary value problems for the equation 
\begin{equation}
\label{eq1}
D_{0y}^{\alpha}u(x,y)+bD_{0y}^{\beta}u(x,y)-u_{xx}(x,y)+cu(x,y)=f(x,y),
\end{equation}
with $0<\alpha=2\beta<1$ was studied in \cite{izvkbnc04} and \cite{daman05}, and in 
the papers \cite{du15-5} and \cite{du15-9} in general case $0<\alpha=2\beta<2$.
For equation (\ref{eq1}) the fundamental solution was constructed, its properties were analyzed, 
a general representation of solutions was obtained, and Cauchy problem was studied in 
\cite{du15-5}, \cite{du15-9}.

For diffusion-wave equation the Green formulas were obtained in papers \cite {Ps-2003} and \cite{pizvran}.
In present paper we obtain the Green formula for operator ${\bf L}$ and construct the solutions of 
basic initial and initial-boundary value problems for equation (\ref{mam1}).

For an extensive bibliography on this subject see \cite{atan2007}, \cite{Huang-2005}, \cite{pizvran}, \cite{MamMon} and references therein.

\section{Statement of problem and main results}\label{sec:2}
\setcounter{section}{2}
\setcounter{equation}{0}\setcounter{theorem}{0}

We denote
$$\Omega=\{(x,y):a_1<x<a_2, 0<y<T\}, 
\Omega_y=\{(t,s):\;a_1<t<a_2,\;0<s<y\}.$$

\begin{definition}\label{Def1}
A function $w=w(x,y;t,s)$ we will call a fundamental solution of equation 
(\ref{mam1}), if it satisfies the following condition.

\noindent
1) For arbitrary fixed point $(x,y)\in \Omega,$ the function $w$ treated as a function 
of the variables $(t,s)$, $s<y,$ satisfies the equation
$${\bf L}^*w(x,y;t,s)\equiv 
\left(D_{ys}^{\alpha} +bD_{ys}^{\beta}-\frac{\partial^2}{\partial t^2}+c\right)w(x,y;t,s)=0; $$
2) The relation 
$$\lim\limits_{s\rightarrow y}\int\limits_{a_1}^{a_2}D_{ys}^{\alpha-1}
w(x,y;t,s)q(t)dt=q(x)$$
holds for each function $q(x)\in C[a_1,a_2].$
\end{definition}

Let us consider the function
\begin{equation} \label{fsfte}
\Gamma(x,y)=\frac{1}{2}\int\limits_{|x|}^{\infty}
h_0(x,\tau)g(y,\tau)d\tau,
\end{equation}
where
$$
h_0(x,\tau)=_0F_1\left[ 1; a(\tau^2-x^2)\right],
\quad
g(y,\tau)=\frac{e^{b_1\tau}}{y}\phi(-\beta,0;-\tau y^{-\beta}),
$$                                                   
$_0F_1(\nu;z)=\sum_{k=0}^{\infty}\frac{\Gamma(\nu)}{\Gamma(k+\nu)}\frac{z^k}{k!}$  
is the confluent hypergeometric limit function \cite [p. 32]{Kilbas-2006},
$\phi(\delta,\mu;z)=\sum_{k=0}^{\infty}\frac{1}{\Gamma(\delta k+\mu)}\frac{z^k}{k!}$  is the Wright function  \cite {Wright-1933}, 
$b_1=-b/2,$ and  $a=b_1^2-c.$

\begin{theorem}\label{Th1}
The function $\Gamma(x-t,y-s)$ is a fundamental solution of equation \eqref{mam1}.
\end{theorem}

Set $J=\{a_1<x<a_2, y=0\};$
here and further
$n\in \{1,2\}$ is a number such that $n-1<\alpha\leq n.$

\begin{definition}\label{Def2}
A regular solution of equation (\ref{mam1}) in a domain $\Omega$ is defined as a  
function $u=u(x,y)$ of the class 
$u\in C(\overline{\Omega}),$ 
$\frac{\partial^{n-1} u}{\partial y^{n-1}}\in C(\Omega\cup J),$ 
$\partial_{0y}^{\alpha}u,\, \partial_{0y}^{\beta}u, 
\,\frac{\partial^2}{\partial x^2}u \in C(\Omega),$
satisfying equation (\ref{mam1}) at all points $(x,y)\in \Omega.$
\end{definition}

\begin{theorem}\label{Th2}
Let $f(x,y)\in C(\overline{\Omega}),$ 
and $\tau_k(x)\in C[a_1,a_2],$ $k=1,n.$
Then any regular solution $u(x,y)$ of equation (\ref{mam1}) in the domain $\Omega$
satisfying the condition
\begin{equation} \label{oos121}
\frac{\partial^{k-1}}{\partial y^{k-1}}u(x,0)=\tau_k(x), 
\quad a_1< x<a_2, \quad k=\overline{1,n}, 
\end{equation}
and such that $u_x\in C(a_1\leq x\leq a_2, 0<y<T),$ $u_x(a_1,y),$ $u_x(a_2,y)\in L[0,T],$
can be represented in the form
$$u(x,y)=\sum\limits_{k=1}^{n}\int\limits_{a_1}^{a_2}
\tau_k(t)\left[D_{0y}^{\alpha-k}+(2-k)bD_{0y}^{\beta-k}\right]G(x,y;t,0)dt+$$
$$+\sum\limits_{i=1}^{2}(-1)^{i}\int\limits_{0}^{y}
[G(x,y;a_i,s)u_t(a_i,s)-G_t(x,y;a_i,s)u(a_i,s)]ds+$$
\begin{equation} \label{kkkkk}
+\int\limits_{a_1}^{a_2}\int\limits_0^y
[G(x,y;t,s)f(t,s)-u(t,s)h(x,y;t,s)]dt ds,
\end{equation}
where
$s^{\alpha-n}h(x,y;t,s)\in L(\Omega_y),$
$G(x,y;t,s)=\Gamma(x-t,y-s)-V(x,y;t,s),$ 
and
$V\equiv V(x,y;t,s)$ 
is a solution of the equation 
\begin{equation} \label{lvh}
{\bf L}^*V(x,y;t,s)=h(x,y;t,s)
\end{equation}
in the class                                                            
$\frac{\partial^2}{\partial t^2}V\in C(\Omega\times\Omega_y)\cup L(\Omega_y),$ 
$V, \frac{\partial}{\partial t}V, D_{ys}^{\alpha}V, D_{ys}^{\beta}V
\in C(\Omega\times\overline{\Omega}_y),$
with the condition
\begin{equation} \label{oos12}
\lim\limits_{s\to y}D_{ys}^{\beta k-i}V(x,y;t,s)=0, \quad
k=1,n, \quad i=1,k.
\end{equation}
\end{theorem}

{\bf Cauchy problem.} 
In the domain $\Omega_0=\{(x,y): x\in \mathbb{R}, 0<y<T\},$ find a solution 
$u(x,y)$ of equation (\ref{mam1}) with the condition
\begin{equation}\label{cpc}
\frac{\partial^{k-1}}{\partial y^{k-1}} u(x,0)=\tau_k (x),
\quad -\infty<x<\infty,
\quad k=1,n,
\end{equation} 
where the $\tau_k(x)$ $k=1,n$ are given functions. 

\begin{definition}\label{Def3}
By $C^{1,q}[a_1, a_2]$ we denote the space of continuously
differentiable functions on $[a_1, a_2]$ which first derivatives 
satisfy the H\"older condition with exponent $q.$
\end{definition}

\begin{theorem}\label{Th3}
Let the functions 
$\tau_k(x)$ $(k=1,n)$ satisfy the conditions
$$\begin{array}{c}
\tau_n(x)\in C(\mathbb{R}); \\
\tau_1(x)\in C^{1,q}(\mathbb{R}),\quad q>\frac{1-\beta}{\beta}, 
\quad \mbox{for }\, n=2;
\end{array} $$
let $f(x,y)\in C(\bar\Omega_0),$
\begin{equation} \label{mam4}
f_y(x,y)\in C(\Omega_0)\cup L(\Omega_0),
\end{equation}
and the relations 
\begin{equation} \label{neq4}
\tau_k (x)=O(\exp(\rho |x|^{\varepsilon})), \quad 
f(x,y)=O(\exp(\rho |x|^{\varepsilon})),  
\end{equation}
where $\varepsilon=\frac{1}{1-\beta}$ and
$\rho<(1-\beta)\left(\frac{\beta}{T}\right)^{\frac{\beta}{1-\beta}},$
hold as $|x|\rightarrow \infty.$  
Then the function
$$ u(x,y)=\sum\limits_{k=1}^{n}\int\limits_{-\infty}^{+\infty}
\tau_k(t)\left[D_{0y}^{\alpha-k}+(2-k)bD_{0y}^{\beta-k}\right]\Gamma(x-t,y)dt+$$
\begin{equation} \label{scpfte}
+\int\limits_0^y \int\limits_{-\infty}^{+\infty}\Gamma(x-t,y-s)f(t,s)dtds
\end{equation}
is a regular solution of the Cauchy problem.
\end{theorem}

\begin{remark}\label{rem1}
For $ 2/3 <\alpha \leq 1$ instead of the condition (\ref{mam4}) of Theorem \ref{Th3} is sufficient 
to require that the function $ f(x, y) $ satisfies a H\"older condition with respect to the 
variable $ x $ with exponent
$q>\frac{1-2\beta}{\beta}.$
\end{remark}

\begin{theorem}\label{Th4}
There exist at most one regular solution $u(x,y)$ of Cauchy problem  in the class of 
functions satisfying the condition
\begin{equation} \label{ndou}
u(x,y)=O(\exp(k |x|^{\varepsilon}))
\quad as \quad |x|\rightarrow \infty
\end{equation} 
for some  $k>0.$ 
\end{theorem}

In the section 6 we will give the statements of the initial boundary value problems for the equation (\ref{mam1}),
and representations of their solutions and of the Green functions.

\section{General representation of the solution }
\label{sec:3}

\setcounter{section}{3}
\setcounter{equation}{0}\setcounter{theorem}{0}

For the proof of Theorems \ref{Th1}, \ref{Th2} and other assertions, some auxiliary statements which were
derived in \cite{du15-5} for the equation (\ref{eq1}) are needed.
Since for function $\Gamma(x,y)$ the relation 
\begin{equation} \label{veq3}
D_{0y}^{\nu}\Gamma(x,y)=\partial_{0y}^{\nu}\Gamma(x,y),
\end{equation}      
holds for all $\nu\in {\mathbb R}$ and for all $(x,y)$ such that $x^2+y^2\not=0,$
then the following lemmas hold also for equation (\ref{mam1}).

\begin{lemma}\label{lem2}
The relations  
\begin{equation} \label{keq24}
D^{\delta}_{0y}D^{\nu}_{0y}\Gamma(x,y)=D^{\nu}_{0y}D^{\delta}_{0y}\Gamma(x,y)
=D^{\delta+\nu}_{0y}\Gamma(x,y)
\end{equation}
hold for arbitrary $\delta,\nu \in \mathbb{R}.$ 
\end{lemma}

\begin{lemma}\label{lem3}
The relation holds 
$$2\frac{\partial^m}{\partial x^m}D_{0y}^{\nu}\Gamma(x,y)=
\int\limits_{|x|}^{\infty}D_{0y}^{\nu-\beta k}g(y,\tau)
\frac{\partial^m}{\partial x^m}L_1^kh_0(x,\tau)d\tau-
$$
$$- \mathop{\rm sgn}(m)\mathop{\rm sgn}(x) \sum\limits_{j=1}^{m}
\frac{\partial^{j-1}}{\partial x^{j-1}}\left[D_{0y}^{\nu-\beta k}g(y,|x|)
\left(\frac{\partial^{m-j}}{\partial x^{m-j}}L_1^kh_0(x,\tau)\right)\big|_{\tau=|x|}
\right]+                                             
$$
$$+\mathop{\rm sgn}(k)\sum\limits_{i=1}^{k}\frac{\partial^m}{\partial x^m}
\left[D_{0y}^{\nu-\beta i}g(y,|x|)
\left(L_1^{i-1}h_0(x,\tau)\right)\big|_{\tau=|x|}\right],
$$
where  $k=0$ for $\nu\leq 0, $ and $k$ such that $\beta (k-1)\nu\leq \beta k$ $(k\in \mathbb{N})$
for $\nu>0,$
is satisfied for arbitrary $m\in \mathbb{N}\cup \{0\}$ and $\nu\in \mathbb{R},$
and $L_1=\frac{\partial}{\partial \tau}+b_1$
is the differential operator of first order.
\end{lemma}

\begin{lemma}\label{lem4}
The estimate 
\begin{equation} \label{kdxdyg}
\left|\frac{\partial^m}{\partial x^m}D_{0y}^{\nu}\Gamma(x,y)\right| 
\leq C |x|^{-\theta}y^{\beta(1-m+\theta)-\nu-1}, 
\quad \theta\geq 0, 
\end{equation}
where $C$ is a positive constant, holds for arbitrary
 $m\in \mathbb{N}\cup \{0\}$ and $\nu\in \mathbb{R}.$ 
\end{lemma}

\begin{remark}\label{rem21}
It follows from Lemma \ref{lem4} that
\begin{equation} \label{ldyg0}
\lim\limits_{y\to 0}D_{0y}^{\mu}\Gamma(x,y)=0, \quad \forall \mu\in {\mathbb R}, \quad x\not=0.
\end{equation}
Hence, by virtue of the next connecting Rimman-Liouville and Caputo derevatives fomula
$$\partial_{ay}^{\nu}g(s)=D_{ay}^{\nu}g(s)-\sum\limits_{k=0}^{n-1}
\frac{g^{(k)}(a)}{\Gamma(k-\nu+1)}|y-a|^{k-\nu},$$
where $n-1<\nu\leq n,$ $n\in \mathbb{N},$  we obtain the equality (\ref{veq3}).
\end{remark}

\begin{lemma}\label{lem5}
The relation ${\bf L}\Gamma(x,y)=0$
holds for all $(x,y)\in\Omega.$
\end{lemma}

\begin{lemma}\label{lem6}
The relation 
${\bf L^*}\Gamma(x-t,y-s)=0 $
holds for all $(t,s)\in\Omega_y$ and for fixed $(x,y)\in \Omega.$ 
\end{lemma}

\begin{lemma}\label{lem8}
The relation
$$\lim\limits_{s\rightarrow y}\int\limits_{x_1}^{x_2}q(t)D^{\alpha-1}_{ys}
\Gamma(x-t,y-s)dt=q(x), \quad x_1<x<x_2.
$$
holds for each function $g(x)\in C[x_1,x_2],$ $a_1\le x_1<x_2\le a_2,$ 
\end{lemma}

\begin{lemma}\label{lem9} 
The relation
$$\lim\limits_{t\to x\pm 0}\int\limits_{\delta}^yp(s)
\frac{\partial}{\partial x}\Gamma(x-t,y-s)ds=\mp \frac{1}{2}p(y), \quad 0\leq\delta< y,
$$
holds for each function $p(y)\in C[\delta,T].$ 
\end{lemma}

Lemmas \ref{lem6} and \ref{lem8}, together with Remark \ref{rem21}, imply Theorem \ref{Th1}.

Let us prove Theorem~\ref{Th2}.

\proof
Let $(x,y)$ is some fixed point in $\Omega,$ 
and $\delta,$ $\varepsilon,$ and $r$ be positive real numbers such that 
$y>\delta>0,$
$a_1+\varepsilon<x-r,$ and $x+r<a_2-\varepsilon.$
By using the relations \cite[p. 34]{n}
$$\int\limits_{\delta}^{y}f(s)D_{\delta s}^{\nu}h(s)ds=
\int\limits_{\delta}^{y}h(s)D_{ys}^{\nu}f(s)ds, \quad 
D_{0s}^{\nu}h(s)=I_{0\delta}^{\nu}h(s)+D_{\delta s}^{\nu}h(s),
$$ 
where $\nu<0,$
$$I_{0\delta}^{\nu}h(s)=
\frac{1}{\Gamma(-\nu)}\int\limits_0^{\delta}h(\xi)(s-\xi)^{-\nu-1}d\xi, 
$$
the integration by parts formula, and by virtue of (\ref{oos12}), 
we transform the integral
$$
\int\limits_{\delta}^{y}G(x,y;t,s)
\left(\partial_{0s}^{\alpha}+b\partial_{0s}^{\beta}\right)u(t,s)ds$$
$$=\int\limits_{\delta}^{y}G(x,y;t,s)
\left[\left(I_{0\delta}^{\alpha-n}+D_{\delta s}^{\alpha-n}\right)
\frac{\partial^n}{\partial s^n}+b\left(I_{0\delta}^{\beta-1}+D_{\delta s}^{\beta-1}\right)
\frac{\partial}{\partial s}\right]u(t,s)ds$$
$$
=\int\limits_{\delta}^{y}u(t,s)\left(D_{ys}^{\alpha}+bD_{ys}^{\beta}\right)
G(x,y;t,s)ds -\sum\limits_{k=1}^{n}\frac{\partial^{k-1}}{\partial s^{k-1}}u(t,s)
D_{ys}^{\alpha-k}G(x,y;t,s)\big |_{s=\delta}$$
$$
-bu(t,s)D_{ys}^{\beta-1}G(x,y;t,s)\big |_{s=\delta}
+\int\limits_{\delta}^{y}G(x,y;t,s)
\left(I_{0\delta}^{\alpha-n}\frac{\partial^n}{\partial s^n}
+bI_{0\delta}^{\beta-1}\frac{\partial}{\partial s}\right)
u(t,s)ds.
$$

Set $B_x^{\varepsilon}=\{t: t\in[a_1+\varepsilon, x-r]\cup [x+r, a_2-\varepsilon]\}.$
By virtue of the last relation and the relation $Gu_{tt}-uG_{tt}=(Gu_t-uG_t)_t,$ we have
$$\int\limits_{\delta}^yds
\int\limits_{B_x^{\varepsilon}}
[G(x,y;t,s){\bf L}u(t,s)-u(t,s){\bf L}^*G(x,y;t,s)]dt$$
$$
=-\int\limits_{B_x^{\varepsilon}}
\left[\sum_{k=1}^{n}\frac{\partial^{k-1}}{\partial s^{k-1}}u(t,s) 
D_{ys}^{\alpha-k}G(x,y;t,s)
+bu(t,s)D_{ys}^{\beta-1}G(x,y;t,s)\right]\big |_{s=\delta}dt
$$
$$+\sum\limits_{i=1}^{2}(-1)^{i+1}\int\limits_{\delta}^y
[G(x,y;\bar a_i,\delta)u_t(\bar a_i,s)-
G_t(x,y;\bar a_i,\delta)u(\bar a_i,s)]ds$$
$$-\int\limits_{\delta}^y[G_t(x,y;x+r,s)u(x+r,s)-G_t(x,y;x-r,s)u(x-r,s)]dt$$
$$+\int\limits_{\delta}^y[G(x,y;x+r,s)u_t(x+r,s)-G(x,y;x-r,s)u_t(x-r,s)]dt$$
\begin{equation} \label{mam51}
+\int\limits_{B_x^{\varepsilon}}
\int\limits_{\delta}^y G(x,y;t,s)
\left[I_{0\delta}^{\alpha-n}\frac{\partial^n}{\partial s^n}+
b I_{0\delta}^{\beta-1}\frac{\partial}{\partial s}\right]u(t,s)ds,
\end{equation}
where $\bar a_i=a_i-(-1)^i\varepsilon,$  $(i=1,2).$

Since $u(x,y), u_y(x,y)\in C((a_1,a_2)\times[0;T]),$ then the estimates
\cite{pizvran}
$$
\left|I_{0\delta}^{\beta-1}\frac{\partial}{\partial s}u(t,s)\right|\leq 
C_0\left[(s-\delta)^{-\beta}+\omega_0(\delta)s^{-\beta}-s^{-\beta}\right],$$
$$
\left|I_{0\delta}^{\alpha-n}\frac{\partial^n}{\partial s^n}u(t,s)\right|\leq 
C_n\left[(s-\delta)^{n-1-\alpha}+\left(\omega_{n-1}(\delta)-1\right)s^{n-1-\alpha}\right],$$
are valid, where 
$\omega_i(\delta)=\sup\limits_{t\in J}\left|\frac{\partial^i}{\partial s^i}u(t,\delta)
-\frac{\partial^i}{\partial s^i}u(t,0)\right|,$ $i=0,n-1.$

Let us compute the following integral
$$I_{\gamma,\nu}(\omega,\delta;y)=\int\limits_{\delta}^{y}
\left[(s-\delta)^{-\gamma}+\omega s^{-\gamma}\right]
(y-s)^{\nu-1}ds, \quad \nu>0, \quad \gamma<1.$$ 
By replacing the integration variables 
$s=y-(y-\delta)\xi,$
we obtain
$$I_{\gamma,\nu}(\omega,\delta;y)=(y-\delta)^{\nu-\gamma}
\int\limits_{0}^{1}(1-\xi)^{-\gamma}\xi^{\nu-1}d\xi$$
$$+\omega y^{-\gamma}(y-\delta)^{\nu}
\int\limits_{0}^{1}\left(1-\frac{y-\delta}{y}\xi\right)^{-\gamma}\xi^{\nu-1}d\xi
=(y-\delta)^{\nu-\gamma}B(1-\gamma,\nu)
$$
\begin{equation} \label{oign}
+\omega y^{-\gamma}(y-\delta)^{\nu}
B(\nu,1)F\left(\gamma,\nu,\nu+1;\frac{y-\delta}{y}\right),
\end{equation}
where
$B(a,b)$ is the beta function,
$F(a,b,c;z)$ is the hipergeometric function.

By using the relation \cite[c. 65]{n}
$$\lim\limits_{z\to 1-0}F\left(\gamma,\nu,\nu+1;z\right)=
\frac{\Gamma(\nu+1)\Gamma(1-\gamma)}{\Gamma(\nu+1-\gamma)},$$
we obtain
\begin{equation} \label{lign}
\lim\limits_{\delta\to 0 \atop \omega\to -1}I_{\gamma,\nu}(\omega,\delta;y)=0.
\end{equation}

By virtue of relation (\ref{oign}) and estimate 
\begin{equation} \label{esg}
\left|\Gamma(x,y)\right|\leq
C |x|^{-\theta}y^{\beta(1+\theta)-1},
\quad \theta\geq 0, 
\end{equation}
which follows from (\ref{kdxdyg}), we see that
$$\left|\int\limits_{\delta}^{y}G(x,y;t,s)
\left(I_{0\delta}^{\alpha-n}\frac{\partial^n}{\partial s^n}
+bI_{0\delta}^{\beta-1}\frac{\partial}{\partial s}\right)
u(t,s)ds\right|
$$
\begin{equation} \label{mamib}
\leq C_0 I_{\beta,\beta(1+\theta)}[\omega_0(\delta)-1,\delta;y]+
C_1 I_{\alpha-n+1,\beta(1+\theta)}[\omega_{n-1}(\delta)-1,\delta;y].
\end{equation}
                                        
By successively letting $r,$ $\delta,$ and $\varepsilon$
tend to zero in relation (\ref{mam51}) and by taking into account the relations
(\ref{oos121}), (\ref{lvh}), (\ref{oign}), (\ref{lign}), (\ref{mamib}), 
relation
$$\lim\limits_{t\rightarrow x+0}\int\limits_{\delta}^{y}u(x,s)\Gamma_t(x-t,y-s)ds
-\lim\limits_{t\rightarrow x-0}\int\limits_{\delta}^{y}u(x,s)\Gamma_t(x-t,y-s)ds=-u(x,y), 
$$
which is valid in view of Lemma~\ref{lem9}
and relation $\lim\limits_{t\to x+0}G=\lim\limits_{t\to x-0}G$ we obtain
(\ref{kkkkk}).

\proofend

\section{Proof of the existence theorem }
\label{sec:4}

\setcounter{section}{4}
\setcounter{equation}{0}\setcounter{theorem}{0}

The rightness of Theorem \ref{Th3} implies from the following 
assertions.

\begin{lemma} \label{lem10}
Let the functions $\tau_k(x)$ $(k=1,n)$  satisfy the conditions
\begin{equation} \label{equstk}
\begin{array}{c}
\tau_n(x)\in C[a_1,a_2]; \\
\tau_1(x)\in C^{1,q}[a_1,a_2],\quad q>\frac{1-\beta}{\beta}, 
\quad \mbox{for }\, n=2.
\end{array}
\end{equation} 
Then function 
\begin{equation} \label{kk5}
u_0(x,y)=\sum\limits_{k=1}^{n}\int\limits_{a_1}^{a_2}
\tau_k(t)\left[D_{0y}^{\alpha-k}+(2-k)bD_{0y}^{\beta-k}\right]\Gamma(x-t,y)dt 
\end{equation}
is a solution of the equation (\ref{mam1}) in the class 
$\partial_{0y}^{\alpha}u_0,$ $\partial_{0y}^{\beta}u_0,$ 
$\frac{\partial^2}{\partial x^2}u_0\in C(\Omega),$ 
and satisfies the conditions  
\begin{equation}\label{neq2220}
\lim\limits_{y \rightarrow 0} 
\frac{\partial^{k-1}}{\partial y^{k-1}} u_0(x,y)=\tau_k (x),
\quad a_1<x<a_2,
\quad k=1,n.
\end{equation}
\end{lemma}

\proof 
It follows from  Lemmas \ref{lem2}, \ref{lem4} and \ref{lem5}
that the function $u_0(x,y)$ is a solution of equation (\ref{mam1}) such that  
$\partial_{0y}^{\alpha}u_0,$ $\partial_{0y}^{\beta}u_0,$ 
$\frac{\partial^2}{\partial x^2}u_0\in C(\Omega ).$

Let us show that  $u_0(x,y)$ satisfies condition (\ref{neq2220}).
The formula \cite{du15-9}
$$
\lim\limits_{y \rightarrow 0}D_{0y}^{k(\beta-1)}\int\limits_{a_1}^{a_2}\tau(t)\Gamma(x-t,s)ds=0, 
$$
holds for all $\tau(x)\in C[a_1,a_2].$
By seting $k=1$ and $k=2$ in the last relation we obtain respectively
\begin{equation} \label{kk6}
\lim\limits_{y \rightarrow 0}\int\limits_{a_1}^{a_2}\tau_2(t)D_{0y}^{\alpha-2}\Gamma(x-t,s)ds=0,
\quad \alpha\in(1,2), 
\end{equation}
\begin{equation} \label{kk7}
\lim\limits_{y \rightarrow 0}\int\limits_{a_1}^{a_2}\tau_1(t)D_{0y}^{\beta-1}\Gamma(x-t,s)ds=0. 
\quad \alpha\in(0,2). 
\end{equation}
By virtue of Lemma \ref{lem8}   
we obtain
\begin{equation} \label{kk8}
\lim\limits_{y \rightarrow 0}
\int\limits_{a_1}^{a_2}\tau_1(t)D_{0y}^{\alpha-1}\Gamma(x-t,s)ds=\tau_1(x), 
\quad \alpha\in(0,2). 
\end{equation}
It follows from relations (\ref{kk7}) and (\ref{kk8}) that for $n=1$
the next relation holds:
$$\lim\limits_{y\to 0}u_0(x,y)=\tau_1(x).$$

Consider the case $n=2.$ 
By virtue of Lemma \ref{lem2} and defferentiation (\ref{kk5}), we obtain
\begin{equation} \label{kk9}
\frac{\partial}{\partial y}u_0(x,y)=\int\limits_{a_1}^{a_2}
\tau_2(t)D_{0y}^{\alpha-1}\Gamma(x-t,y)dt
+\int\limits_{a_1}^{a_2}
\tau_1(t)\left[D_{0y}^{\alpha}+bD_{0y}^{\beta}\right]\Gamma(x-t,y)dt.
\end{equation}

In paper \cite{du15-9} was shoun that the relation
\begin{equation} \label{kk19}
\lim\limits_{y \rightarrow 0}\int\limits_{a_1}^{a_2}
\tau_1(t)\left[D_{0y}^{\alpha}+bD_{0y}^{\beta}\right]\Gamma(x-t,y)dt=0,
\quad \alpha\in(1,2), 
\end{equation}
holds for each function 
$\tau_1(x)\in C^{1,q}[a_1,a_2],$ $q>\frac{1-\beta}{\beta}.$
From Lemma \ref{lem8} we obtain the relation
\begin{equation} \label{kk10}
\lim\limits_{y \rightarrow 0}\int\limits_{a_1}^{a_2}\tau_2(t)D_{0y}^{\alpha-1}\Gamma(x-t,y)dt
=\tau_2(x).
\end{equation}
Thus, the relations (\ref{kk6}) -- (\ref{kk8}) imply that the relation (\ref{neq2220}) holds for $k=1,$ 
and the relations (\ref{kk8}) -- (\ref{kk10}) imply that the relation (\ref{neq2220}) holds for $k=2.$
\proofend

\begin{lemma} \label{lem11}
Let 
$f(x,y)\in C(\bar\Omega),$
$f_y(x,y)\in C(\Omega)\cup L(\Omega).$
Then the function 
\begin{equation}\label{uf}
u_f(x,y)=\int\limits_0^y \int\limits_{a_1}^{a_2}\Gamma(x-t,y-s)f(t,s)dtds 
\end{equation}
is a regular solution of equation (\ref{mam1}) in domain $\Omega$ satisfying 
the homogenous condition (\ref{neq2220}).
\end{lemma}

\proof
For $\alpha\in (0,1]$ we have 
$$\frac{\partial}{\partial y}u_f(x,y)=\frac{\partial}{\partial y}
\int\limits_0^y\int\limits_{a_1}^{a_2}\Gamma(x-t,s)f(t,y-s)dtds$$
$$=\int\limits_{a_1}^{a_2}\Gamma(x-t,y)f(t,0)dt+
\int\limits_0^y\int\limits_{a_1}^{a_2}\Gamma(x-t,y-s)f_s(t,s)dtds.
$$
Thus
$$\left(\partial_{0y}^{\alpha}+b\partial_{0y}^{\beta}\right)u_f(x,y)=
\int\limits_{a_1}^{a_2}\left(D_{0y}^{\alpha-1}+bD_{0y}^{\beta-1}\right)\Gamma(x-t,y)f(t,0)dt 
$$
\begin{equation}\label{eq32}
+\int\limits_0^y\int\limits_{a_1}^{a_2}
\left(D_{ys}^{\alpha-1}+bD_{ys}^{\beta-1}\right)\Gamma(x-t,y-s)f_s(t,s)dtds. 
\end{equation}

For $\alpha\in (1,2),$ by taking into account the relation (\ref{ldyg0}) with $\mu=0,$ we have
\begin{equation}\label{duf}
\frac{\partial}{\partial y}u_f(x,y)=
\int\limits_0^y\int\limits_{a_1}^{a_2}\Gamma_y(x-t,y-s)f(t,s)dtds,
\end{equation}
$$\frac{\partial^2}{\partial y^2}u_f(x,y)=
\int\limits_{a_1}^{a_2}\Gamma_y(x-t,y)f(t,0)dt+
\int\limits_0^y\int\limits_{a_1}^{a_2}\Gamma_y(x-t,y-s)f_s(t,s)dtds.
$$
Thus
$$\left(\partial_{0y}^{\alpha}+b\partial_{0y}^{\beta}\right)u_f(x,y)=
D_{0y}^{\alpha-2}\int\limits_{a_1}^{a_2}\frac{\partial}{\partial y}\Gamma(x-t,y)f(t,0)dt
$$
\begin{equation}\label{eq33}
+D_{0y}^{\alpha-2}\int\limits_0^y\int\limits_{a_1}^{a_2}\Gamma_y(x-t,y-s)f_s(t,s)dtds+ 
D_{0y}^{\beta-1}\int\limits_0^y\int\limits_{a_1}^{a_2}
\Gamma_y(x-t,y-s)f(t,s)dtds.
\end{equation}
The estimate (\ref{kdxdyg})
implies that $\partial^{\alpha}_{0y}u_f, \partial^{\beta}_{0y}u_f\in C(\Omega).$

Let us find the derivative $\frac{\partial^2}{\partial x^2}u_f(x,y).$ 
For this end, we consider the function 
$$F_{\delta_1,\delta_2}(x,y)=\int\limits_0^yds
\int\limits_{B_x}\Gamma_x(x-t,y-s)f(t,s)dt,$$
where $B_x=\{t: t\in[a_1, x-\delta_1]\cup [x+\delta_2, a_2]\}.$
Obviously, that 
$$\lim\limits_{\delta_1\rightarrow 0 \atop \delta_2\rightarrow 0}F_{\delta_1,\delta_2}(x,y)=
\frac{\partial}{\partial x}u_f(x,y)\in C(\Omega).$$
By using Lemma \ref{lem6} we derive the derivative of the function 
$F_{\delta_1,\delta_2}(x,y)$ by variable $x:$ 
$$\frac{\partial}{\partial x}F_{\delta_1,\delta_2}(x,y)=
\int\limits_0^y \Gamma_x(\delta_1,y-s)f(x-\delta_1,s)ds-
\int\limits_0^y \Gamma_x(-\delta_2,y-s)f(x+\delta_2,s)ds$$
$$+\int\limits_0^yds
\int\limits_{B_x}
f(t,s)\left(D_{ys}^{\alpha}+bD_{ys}^{\beta}\right)\Gamma(x-t,y-s)dt$$
\begin{equation}\label{eq34}
+c\int\limits_0^yds
\int\limits_{B_x}
f(t,s)\Gamma(x-t,y-s)dt.
\end{equation}

We denote the third summand in right hand side of (\ref{eq34}) via  $J_{\delta_1,\delta_2}(x,y).$
Then, by integration by parts, for $\alpha\in(0,1]$ $(n=1)$  we have
$$J_{\delta_1,\delta_2}(x,y)=-\int\limits_0^yds
\int\limits_{B_x}
f(t,s)\frac{\partial}{\partial s}
\left(D_{ys}^{\alpha-1}+bD_{ys}^{\beta-1}\right)\Gamma(x-t,y-s)dt$$
$$=\left(D_{0y}^{\alpha-1}+bD_{0y}^{\beta-1}\right)\int\limits_{B_x} 
\left[f(t,0)\Gamma(x-t,y)+
\int\limits_0^y 
f_s(t,s)\Gamma(x-t,y-s)ds\right]dt.$$
From the last relation and (\ref{eq34}) we obtain
$$\frac{\partial^2}{\partial x^2}u_f(x,y)=
\lim\limits_{\delta_1\rightarrow 0 \atop \delta_2\rightarrow 0}
\frac{\partial}{\partial x}F_{\delta_1,\delta_2}(x,y)= -f(x,y)+cu_f(x,y)$$
\begin{equation} \label{dxf1}
+\left(D_{0y}^{\alpha-1}+bD_{0y}^{\beta-1}\right)\int\limits_{a_1}^{a_2} 
\left[f(t,0)\Gamma(x-t,y)+
\int\limits_0^y f_s(t,s)\Gamma(x-t,y-s)ds\right]dt.
\end{equation}

For $\alpha\in(1,2)$ $(n=2)$  we have
$$J_{\delta_1,\delta_2}(x,y)=
\int\limits_0^yds \int\limits_{B_x} f(t,s)
\left(\frac{\partial^2}{\partial s^2}D_{ys}^{\alpha-2}-
b\frac{\partial}{\partial s}D_{ys}^{\beta-1}\right)\Gamma(x-t,y-s)dt$$
$$=D_{0y}^{\alpha-2}\int\limits_{B_x} 
\left[f(t,0)\Gamma_y(x-t,y)+ \int\limits_0^y 
f_s(t,s)\Gamma_y(x-t,y-s)ds\right]dt$$
$$+bD_{0y}^{\beta-1}\int\limits_0^yds \int\limits_{B_x} f(t,s)\Gamma_y(x-t,y-s)dt.
$$
Hence, 
$$\frac{\partial^2}{\partial x^2}u_f(x,y)=
\lim\limits_{\delta_1\rightarrow 0 \atop \delta_2\rightarrow 0}
\frac{\partial}{\partial x}F_{\delta_1,\delta_2}(x,y)=-f(x,y)+cu_f(x,y)$$
$$+D_{0y}^{\alpha-2}\int\limits_{a_1}^{a_2}f(t,0)\Gamma_y(x-t,y)dt 
+ D_{0y}^{\alpha-2}\int\limits_0^y \int\limits_{a_1}^{a_2}
f_s(t,s)\Gamma_y(x-t,y-s)dtds$$
\begin{equation} \label{dxf2}
+bD_{0y}^{\beta-1}\int\limits_0^yds \int\limits_{a_1}^{a_2} f(t,s)\Gamma_y(x-t,y-s)dt.
\end{equation}

It follows from the relations (\ref{uf}), (\ref{eq32}), (\ref{eq33}), (\ref{dxf2}) and (\ref{kdxdyg}) that
$${\bf L}u_f(x,y)=f(x,y),$$
moreover 
$\frac{\partial^2}{\partial x^2}u_f,$ $\partial_{0y}^{\alpha}u_f,$ 
$\partial_{0y}^{\beta}u_f \in C(\Omega).$

From estimate (\ref{esg}) and the relation (\ref{uf})  it is easy to see that 
$u_f(x,0)=0,$ for $\alpha\in (0,2).$
By virtue of the relation (\ref{duf}) and the estimate 
$$\left|\frac{\partial}{\partial y}\Gamma(x,y)\right|
\leq C |x|^{-\theta}y^{\beta+\beta\theta-2}, \quad \theta\geq 0, 
$$
which follows from (\ref{kdxdyg}), and by taking into account the fact that 
the inequality  $\frac{1-\beta}{\beta}<1$ holds for $\beta>\frac{1}{2},$  we obtain that the relation 
$\frac{\partial}{\partial y}u_f(x,0)=0$ holds for $\alpha\in (1,2).$
\proofend

The following two lemmas are proved in paper \cite{du15-9}.

\begin{lemma} \label{lem12}
The asymptotic relations 
\begin{equation} \label{neq19}
\frac{\partial^m}{\partial x^m}D_{0y}^{\nu}\Gamma(x,y)=
O\left(\exp\left(-\rho_{\beta}(y)|x|^{\varepsilon}\right)\right),
\quad \beta-\nu-\beta m\geq 0, 
\end{equation}
$$\frac{\partial^m}{\partial x^m}D_{0y}^{\nu}\Gamma(x,y)=
O\left(\exp\left(-\rho^*_{\beta}(y)|x|^{\varepsilon}\right)\right), 
\quad \beta-\nu-\beta m< 1, 
$$
hold as $|x|\rightarrow \infty,$ 
where 
$m\in {\mathbb N}\cup \{0\},$  $\nu\in {\mathbb R},$
$\rho_{\beta}(y)=(1-\beta)\left(\beta y^{-1}\right)^{\frac{\beta}{1-\beta}},$
$\varepsilon=\frac{1}{1-\beta},$ 
$\rho^*_{\beta}(y)=\rho_{\beta}(y)\theta-\delta,$ 
and $\theta$ and $\delta$ are arbitrary numbers in the interval $]0,1[.$
\end{lemma}

\begin{lemma} \label{lem13}
Let  
$\tau_k (x) \in C(-\infty,+\infty),$   $y^{n-\alpha}f(x,y)\in C(\bar\Omega),$
and let the relations  (\ref{neq4})
hold as $|x|\rightarrow \infty.$  
Then the following inclusions hold
$$
\int\limits_{x\pm\delta}^{\pm\infty}\frac{\partial^m}{\partial x^m}D_{0y}^{\nu}
\Gamma(x-t,y)\tau_k(t)dt,
\int\limits_0^y\int\limits_{x\pm\delta}^{\pm\infty}
\frac{\partial^m}{\partial x^m}D_{ys}^{\nu}
\Gamma(x-t,y-s)f(t,s)dtds\in C(\Omega),$$
where $m\in {\mathbb N}\cup \{0\}$ and  $\nu\in {\mathbb R}.$    
\end{lemma}

\begin{remark} \label{Rem2}
If the additional conditions (\ref{neq4}) are satisfied,
then by using Lemmas \ref{lem12}, \ref{lem13} and Remark \ref{rem21}, one can readilly see that Lemmas \ref{lem10} and \ref{lem11} hold
for the case in which both $a_1$ and $a_2$ or one of them is infinite.
\end{remark}

\section{Proof of the uniqueness theorem }
\label{sec:5}

\setcounter{section}{5}
\setcounter{equation}{0}\setcounter{theorem}{0}

Let us prove Theorem~\ref{Th4}.

\proof
Let  $h_r(t)$ be a function such that $h_r(t)\in C^2(\mathbb R),$ 
$$ h_r(t)=
\left \{
\begin{array}{ll}
1, 			& |t|\leq r,\\
0,			& |t|\geq r+1.
\end{array}
\right.$$
$0\leq h_r(t)\leq 1,$ $|h'_r(t)|\leq C,$ and $|h''_r(t)|\leq C,$ 
where $C={\rm const}>0.$

It follows from Lemma \ref{lem6} that 
\begin{equation} \label{nlhg}
{\bf L^*}h_r(t)\Gamma(x-t;y-s)=-h''_r(t)\Gamma(x-t,y-s)-2h'_r(t)\Gamma_t(x-t,y-s).
\end{equation} 
By virtue of the relation (\ref{nlhg}) and Theorem \ref{Th2},  
we find that the regular solution of the homogeneous problem (\ref{mam1}), (\ref{cpc}) in the domain
$$\Omega_r=\{(x,y):|x|<r+1, 0<y<T\}$$
can be represented in the form 
\begin{equation} \label{nsol}
u(x,y)=
\int\limits_0^y\left(\int\limits_{-r-1}^{-r}+\int\limits_{r}^{r+1}\right)
[h''_r\Gamma+2h'_r\Gamma_t]u(t,s)dtds.
\end{equation}
By taking into account the estimate (\ref{neq19}) and the fact that, by virtue of 
(\ref{ndou}), the inequality 
$$|u(x,y)|\leq C \exp(k |x|^{\varepsilon}),$$
which holds for sufficiently large $r,$ from relation (\ref{nsol}) we obtain 
$$|u(x,y)|\leq C
\left(\int\limits_{-r-1}^{-r}+\int\limits_{r}^{r+1}\right)\int\limits_0^y
\exp\left[-\rho_{\beta}(s)|x-t|^{\varepsilon}+k|t|^{\varepsilon}\right]dsdt
$$
\begin{equation} \label {ipm1}
=I_+(x,y)+I_-(x,y),
\end{equation}
where
$\rho_{\beta}(y)=\sigma y^{-\beta\varepsilon},$
$$I_{\pm}(x,y)=C\int\limits_{r}^{r+1}\int\limits_0^y
\exp[-\rho_{\beta}(s)|x\pm t|^{\varepsilon}+k|t|^{\varepsilon}]dsdt.$$
It is easy to see that, for $0<t<\infty,$ the function 
$$f_{\pm}(t)=-\rho_1|x\pm t|^{\varepsilon}+k|t|^{\varepsilon}$$
has the unique maximum
\begin{equation} \label{max55}
\max f_{\pm}(t)=\mu |x|^{\varepsilon},
\quad
\mu=k\rho_1\left(\rho_1^{\gamma}-k^{\gamma}\right)^{-\frac{1}{\gamma}},
\quad
\gamma=\frac{1}{\varepsilon-1}.
\end{equation}

The inequalities 
$k<\rho_1=\sigma T_1^{-\beta\varepsilon}<\rho_{\beta}(y)=\sigma y^{-\beta\varepsilon},$
hold for all
$y<T_1=\beta[(1-\beta)/\sigma]^{(1-\beta)/\beta},$
therefore, by (\ref{max55}), we have
$$I_{\pm}(x,y)=C\int\limits_{r}^{r+1}\int\limits_0^y
\exp\left[-\rho_2(s)|x\pm t|^{\varepsilon}\right]
\exp\left[-\rho_1|x\pm t|^{\varepsilon}+k|t|^{\varepsilon}\right]dsdt
$$
$$\leq Ce^{\mu |x|^{\varepsilon}}
\int\limits_0^yds\int\limits_{r}^{r+1}
\exp\left[-\rho_2(s)|x\pm t|^{\varepsilon}\right]dt$$
\begin{equation} \label{ipm2}
\leq Ce^{\mu |x|^{\varepsilon}}\int\limits_0^y
[F(r\pm x,s)-F(r+1\pm x,s)]ds,
\end{equation}
where $\rho_2(s)=\rho_{\beta}(s)-\rho_1,$ and
$F(x,s)=\frac{\rho_2(s)}{\varepsilon}x^{1-\varepsilon}\exp[-\rho_2(s)x^{\varepsilon}].$
Since, for fixed $x,$ the maximum of the function $F(x,s)$ is attained at
$\rho_2(s)=x^{-\varepsilon},$ i.e.,
$$\max\limits_{0<s<T_0} F(x,s)=\frac{1}{\varepsilon}x^{1-2\varepsilon}\exp(-1),$$
it follows from (\ref{ipm2}) that
$$0\leq I_{\pm}(x,y)\leq\frac{C}{\varepsilon}e^{\mu |x|^{\varepsilon}-1}
y\left[(r\pm x)^{1-2\varepsilon}+(r+1\pm x)^{1-2\varepsilon}\right].
$$
This, together with the relations
$1-2\varepsilon=-\frac{1+\beta}{1-\beta}<0$
implies that 
\begin{equation} \label {ipm3}
\lim\limits_{r\rightarrow \infty} I_{\pm}(x,y)=0.
\end{equation}
Now it follows from the estimate (\ref{ipm1}) and the relation (\ref{ipm3}) that
$u(x,y)\equiv 0$ in the domain
$$\Omega_1=\{(x,y):x\in \mathbb{R},\, 0<y\leq T_1\}.$$

Let us show that $u(x,y)\equiv 0$ in the strip 
$$\Omega_2=\{(x,y):x\in \mathbb{R},\, T_1<y\leq 2T_1\}.$$
Consider the function $v(x,z)=u(x,T_1+z),$ where $z=y-T_1.$
Since $u(x,y)\equiv 0$ in $\Omega_1,$ we have
$$D_{0y}^{\nu}u(x,y)=D_{T_1y}^{\nu}u(x,y)=D_{0z}^{\nu}v(x,z).$$
Hence it follows that for arbitrary 
$(x,z)\in \Omega_1,$
the function $v(x,z)$ is a solution of the homogeneous Cauchy problem.
Therefore, the above argument implies that $v(x,z)\equiv 0$ for all
$(x,z)\in \Omega_1.$ 
Consequently, $u(x,y)\equiv 0$  for all $(x,y)\in \Omega_2.$ 

In the same way, one can show that $u(x,y)\equiv 0$ in
$$\Omega_3=\{(x,y):x\in \mathbb{R},\, 2T_1<y\leq 3T_1\}$$
and so on.
\proofend

\section{Green functions and representations of solutions }
\label{sec:6}

\setcounter{section}{6}
\setcounter{equation}{0}\setcounter{theorem}{0}

In this section we give the statements of the basic boundary value problems for equation (\ref{mam1})
and present their solutions in term of the appropriate Green functions.

{\bf Problem ${\bf P_i (i=0,1).}$}
In the domain $\Omega_+=\{(x,y): x>a_1, 0<y<T\},$ find a solution 
$u(x,y)$ of equation (\ref{mam1}) with the conditions 
$$
\frac{\partial^{k-1}}{\partial y^{k-1}} u(x,0)=\tau_k (x),
\quad a_1<x<\infty,
\quad k=1,n,
$$
$$\frac{\partial^i}{\partial x^i}u(a_1,y)=\varphi(y), \quad  0<y<T,
$$
where the $\tau_k(x)$ $(k=1,n),$ and $\varphi(y)$ are given functions. 

{\bf Problem ${\bf P_{ij} (i,j=0,1).}$}
In the domain $\Omega=\{(x,y): a_1<x<a_2, 0<y<T\},$ find a solution 
$u(x,y)$ of equation (\ref{mam1}) with the conditions (\ref{oos121}) and
$$\frac{\partial^i}{\partial x^i}u(a_1,y)=\varphi_1(y), \quad  
\frac{\partial^j}{\partial x^j}u(a_2,y)=\varphi_2(y), \quad  0<y<T,
$$
where the $\tau_k(x)$ $k=1,n,$ $\varphi_1(y)$ and $\varphi_2(y)$ are given functions. 

Thus the problems $P_0$ and $P_1$ are the boundary value problems in half-strip with
first and second kind boundary conditions respetively, 
and problems $P_{00},$ $P_{11},$ and $P_{01}$ and $P_{10}$  are the first, second and mixed
boundary value problems for equation (\ref{mam1}) respectively.

\begin{definition}\label{Def4}
The function $G(x,y;t,s)=\Gamma(x-t,y-s)-V(x,y;t,s),$ which satisfied the conditions
$$\lim\limits_{t \rightarrow a_1}\frac{\partial^i}{\partial t^i}G(x,y;t,s)=0, \quad y\not=s;
\quad
\lim\limits_{t \to \infty}\frac{\partial^{k}}{\partial t^{k}}G(x,y;t,s)=0, \quad k=0,1,
$$
or
$$
\lim\limits_{t \rightarrow a_1}\frac{\partial^i}{\partial t^i}G(x,y;t,s)=0,\ 
\lim\limits_{t \rightarrow a_2}\frac{\partial^j}{\partial t^j}G(x,y;t,s)=0, \ y\not=s,
$$
where  $V\equiv V(x,y;t,s)$ is a function  under condition of Theorem \ref{Th2},
with $h(x,y;t,s)\equiv 0,$
we will call {\it a Green function of the problem} $P_i$ 
or $P_{ij}$ for equation (\ref{mam1}) respectively, $(i,j=0,1).$
\end{definition}

Lemmas \ref{lem4} and \ref{lem6} imply the following assertion.

\begin{theorem}\label{Th5}
The function 
$$G_{i}(x,y;t,s)=\Gamma(x-t,y-s)+(-1)^{i+1}\Gamma(x+t-2a_1,y-s),$$
is a Green function of problem $P_{i},$ $(i=0,1).$ 
\end{theorem}

\begin{theorem}\label{Th6}
The function 
$$G_{ij}(x,y;t,s)=\sum\limits_{m=-\infty}^{+\infty}
(-1)^{(i+j)m}\bigg[\Gamma(X_1^m,y-s)+(-1)^{i+1}\Gamma(X_2^m,y-s)\bigg],$$
where 
$$X_1^m=2m(a_2-a_1)+x-t, \quad X_2^m=2m(a_2-a_1)+x+t-2a_1,$$
is a Green function of problem $P_{ij},$ $(i,j=0,1).$ 
\end{theorem}

We denote via $u_{i}(x,y)$ and $u_{ij}(x,y)$ the solutions of problems $P_{i}$ and $P_{ij}$ respectively,
$(i,j=0,1).$  

It follows from theorems \ref{Th2}, \ref{Th5} and \ref{Th6} that the function 
$u_i(x,y)$ has a form
$$u_{i}(x,y)=\sum\limits_{k=1}^{n}\int\limits_{a_1}^{\infty}
\tau_k(t)\left[D_{0y}^{\alpha-k}+(2-k)bD_{0y}^{\beta-k}\right]G_{i}(x,y;t,0)dt$$
\begin{equation} \label{ui}
+(-1)^{i}\int\limits_{0}^{y}
\frac{\partial^{1-i}}{\partial t^{1-i}}G_{i}(x,y;a_1,s)\varphi_1(s)ds+
\int\limits_{0}^y\int\limits_{a_1}^{\infty} G_{i}(x,y;t,s)f(t,s)dtds,
\end{equation}
and function $u_{ij}(x,y)$ has a form
$$u_{ij}(x,y)=\sum\limits_{k=1}^{n}\int\limits_{a_1}^{a_2}
\tau_k(t)\left[D_{0y}^{\alpha-k}+(2-k)bD_{0y}^{\beta-k}\right]G_{ij}(x,y;t,0)dt $$
$$+\int\limits_{0}^{y}
\left[(-1)^{i}\frac{\partial^{1-i}}{\partial t^{1-i}}G_{ij}(x,y;a_1,s)\varphi_1(s)-
(-1)^{j}\frac{\partial^{1-j}}{\partial t^{1-j}}G_{ij}(x,y;a_2,s)\varphi_2(s)\right]ds$$
$$+
\int\limits_{0}^y\int\limits_{a_1}^{a_2} G_{ij}(x,y;t,s)f(t,s)dtds.$$

The representation (\ref{ui}) holds if the functions $\tau_k(x)$ $(k=1,n)$ and $f(x,y)$
satisfy the condition (\ref{neq4}) as $x\to +\infty.$

\begin{remark}\label{rem61}
For $b=c=0$ the equation (\ref{mam1}) reduces to diffusion-wave equation
$$\partial_{0y}^{\alpha}u(x,y)-\frac{\partial^2}{\partial x^2}u(x,y)=f(x,y).$$ 
In this case the function (\ref{fsfte}) takes the form 
$$\Gamma(x,y)=\frac{1}{2y}\int\limits_{|x|}^{\infty}\phi(-\beta,0;-\tau y^{-\beta})=\frac{y^{\beta-1}}{2}\phi(-\beta,\beta;-|x|y^{-\beta})$$
and the representation (\ref{scpfte})  of Cauchy problem solution takes the form
$$u(x,y)=\frac{1}{2}\sum\limits_{k=1}^{n}\int\limits_{-\infty}^{+\infty}\tau_k(t)D_{0y}^{\alpha-k}y^{\beta-1}\phi(-\beta,\beta;-|x-t|y^{-\beta})dt$$
$$+\frac{1}{2}\int\limits_{-\infty}^{+\infty}dt\int\limits_{0}^{y}f(t,s)(y-s)^{\beta-1}\phi(-\beta,\beta;-|x-t|(y-s)^{-\beta})ds.$$
These representations coincides with the representations which were obtained in works \cite{PsMon} and \cite{pizvran}.

For $f(x,y)\equiv 0,$ $\tau_2(x)\equiv 0$ the last relation can be expressed as
$$u(x,y)=\frac{1}{2}\int\limits_{-\infty}^{+\infty}\tau_1(t)D_{0y}^{\alpha-1}y^{\beta-1}\phi(-\beta,\beta;-|x-t|y^{-\beta})dt$$
$$=\frac{1}{2}\int\limits_{-\infty}^{+\infty}\tau_1(t)y^{-\beta}\phi(-\beta,1-\beta;-|x-t|y^{-\beta})dt.$$
That is consistent with the result obtained in paper \cite{Mainardi-1996}.
\end{remark}

\section{Conclusion }
\label{sec:7}

\setcounter{section}{7}
\setcounter{equation}{0}\setcounter{theorem}{0}

In paper the time-fractional telegraph equation is investigated by the Green function method.
The integral representation of a regular solutions is obtained in the form of the Green formula
for operator generating the equation (\ref{mam1}).
The solutions of the basic initial and initial-boundary value problems are constructed in
terms of fundamental solution and appropriate Green functions.

Comparing the results of this paper with the results of the papers \cite{du15-5} and \cite{du15-9}, we conclude
that for the equations (\ref{mam1}) and (\ref{eq1}) the fundamental solutions, and corresponding
Green functions coincide, but the form of initial conditions, the solutions representations, 
and the classes of regular solutions of corresponding problems is different. 
This is due to the fact that the operators  $\partial_{0y}^{\nu}$ and $D_{0y}^{\nu}$ are equal 
for a certain class of functions,
but in general have a different domains of determination.

All results can be easily extended to the fractional advection-dispersion equation (\ref{fade}).

The obtained representations can be used to solve the non-local boundary value problems for the equation (\ref{mam1}),
and for the realization of models described by mixed type equations, which in one of the parts of the mixed domain, coincide
with equation (\ref{mam1}).

\section*{Acknowledgements}

This work was supported by the Division of Nanotechnologies
and Information Technologies of the
Russian Academy of Sciences, project no. 5 "Fundamental
Problems and Technologies of Epitaxial
Nanostructures and Devices Based on Them."




 \bigskip \smallskip

 \it

\noindent
Institute of Applied Mathematics and Avtomatition \\
"Shortanov" Str., 89 A, \\
Nal'chik -- 360000, RUSSIA  \\[4pt]
e-mail: mamchuev@rambler.ru

\end{document}